\documentclass{amsart}
\usepackage{amssymb, amsmath}
\usepackage[latin1]{inputenc}
\usepackage{hyperref}

\usepackage{color}

\newtheorem{theorem}{Theorem}
\newtheorem{lemma}[theorem]{Lemma}

\newtheorem{remark}[theorem]{Remark}

\newtheorem{definition}[theorem]{Definition}

\begin{document}

\title{Locally conformally flat Lorentzian quasi-Einstein manifolds}
\author{M. Brozos-V\'{a}zquez $\,$ E. Garc\'{i}a-R\'{i}o $\,$ S. Gavino-Fern\'{a}ndez}
\address{MBV: Department of Mathematics, University of A Coru\~na, Spain}
\email{miguel.brozos.vazquez@udc.es}
\address{EGR-SGF: Faculty of Mathematics,
University of Santiago de Compostela,
15782 Santiago de Compostela, Spain}
\email{eduardo.garcia.rio@usc.es $\,\,$ sandra.gavino@usc.es}
\thanks{
M. B.-V. and E. G.-R. are supported by projects MTM2009-07756 and INCITE09 207 151 PR (Spain). S. G.-F. is supported by project MTM2009-14464-C02-01 (Spain)}

\begin{abstract}
We show that locally conformally flat quasi-Einstein manifolds are globally conformally equivalent to a space form or locally isometric to a $pp$-wave or a warped product.
\end{abstract}

\maketitle

\section{Introduction}\label{se:1}

Einstein metrics are important both in mathematics and physics.
Some generalizations like conformally Einstein metrics and Ricci solitons play an important role in understanding the AdS/CFT correspondence and the Ricci flow, respectively \cite{A}, \cite{Cao}. In this paper we focus on another generalization which is related to the study of smooth metric measure spaces. A smooth metric measure space is a Riemannian manifold with a measure which is conformal to the Riemannian one. Formally, it is a triple $(M,g,e^{-f}dvol_g)$, where $M$ is a complete $n$-dimensional smooth manifold with Riemannian metric $g$, $f$ is a smooth real valued function on $M$, and $dvol_g$ is the Riemannian volume density on $(M,g)$. This is also sometimes called a manifold with density.

A natural extension of the Ricci tensor to smooth metric measure spaces is the $m$-Bakry-Emery Ricci tensor
\[
\rho_f^m=\rho+\text{Hes}_f-\frac 1m df\otimes df,\quad \text{for}\quad 0<m\leq \infty.
\]
When $f$ is constant, this is the usual Ricci tensor. We call a quadruple $(M,g,f,-\frac{1}{m})$ quasi-Einstein if it satisfies the equation $\rho_f^m=\lambda g$, for some $\lambda \in \mathbb{R}$. Quasi-Einstein manifolds are a generalization of Einstein metrics and they contain the gradient Ricci solitons as a limit case when $m=\infty$. Moreover they are also closely related to the construction of warped product Einstein metrics. Recall here that a complete classification of warped product Einstein metrics is still an open problem in spite of their interest in describing standard stationary metrics.
Quasi-Einstein metrics and, more generally, the $m$-Bakry-Emery Ricci tensor allows some generalizations of the celebrated Hawking-Penrose singularity theorems and the Lorentzian splitting theorem \cite{case}.

We are going to work in a slightly more general setting. Henceforth let $(M,g)$ be a Lorentzian manifold of dimension $n+2$, with $n\geq 1$. Let $f$ be a smooth function on $M$ and let $\mu\in\mathbb{R}$ be an arbitrary constant. We say that $(M,g,f,\mu)$ is \emph{quasi-Einstein} if there exists $\lambda\in \mathbb{R}$ so that
\begin{equation}\label{q-einstein}
\rho+\text{Hes}_f-\mu df\otimes df=\lambda g.
\end{equation}
If we particularize the elements of this quadruple we get some important families of manifolds. Thus, for example, if $\mu=0$ then quasi-Einstein manifolds correspond to gradient Ricci solitons and if $f$ is constant \eqref{q-einstein} reduces to the Einstein equation.
There are other interesting relations between quasi-Einstein manifolds and some well-known structures:
\begin{enumerate}
\item[(a)] If $(M,g,f,-\frac{1}{n})$ is quasi-Einstein, then the conformal metric $\widetilde{g}=e^{-\frac{2}{n}f}g$ is Einstein since, from the expression of the Ricci tensor of a conformal metric, we get \cite{sakai}
\begin{align*}
\rho_{\widetilde{g}} &=\rho_g+\text{Hes}_f+\frac{1}{n}df\otimes df+\frac{1}{n}(\Delta f-\|\nabla f\|^2)g\\
                     &=\frac{1}{n}(\Delta f-\|\nabla f\|^2+n\lambda)e^{\frac{2}{n}f}\widetilde{g},
\end{align*}
where $\nabla f$ and $\Delta f$ denote the gradient and the Laplacian of $f$. In particular, if $g$ is also locally conformally flat, then $\widetilde{g}$ has constant curvature.

\item[(b)] Let $M\times_f F$ be an Einstein warped product. Then it follows from the expressions of the Ricci tensor of the warped product (see, for example, \cite{FGKU}) that $(M,g,-(\operatorname{dim}\, F)\,\operatorname{log}(f),\frac{1}{\operatorname{dim}\, F})$ is quasi-Einstein.
\end{enumerate}

Locally conformally flat complete quasi-Einstein Riemannian manifolds were recently classified in \cite{CMMR}. The purpose of this work is to describe the local structure of locally conformally flat quasi-Einstein manifolds in the Lorentzian setting. The main result is summarized as follows

\begin{theorem}\label{mainth}
Let $(M,g,f,\mu)$ be a locally conformally flat Lorentzian quasi-Einstein manifold.
\begin{enumerate}

\item [(i)] If $\mu=-\frac{1}{n}$, then $(M,g)$ is globally conformally equivalent to a space form.

\item [(ii)] If $\mu\neq -\frac{1}{n}$, then

    \begin{enumerate}
        \item In a neighborhood of any point where $\|\nabla f\|\neq 0$, $M$ is locally isometric to a warped product $I\times_{\phi} F$, where $I$ is a real interval and $F$ is a $(n+1)$-dimensional fiber of constant sectional curvature.

        \item If $\|\nabla f\|=0$, then $(M,g)$ is locally isometric to a plane wave, i.e., $(M,g)$ is locally isometric to $\mathbb{R}^2\times\mathbb{R}^n$ with metric
\[
            g=2dudv +H(u,x_1,\dots,x_n) du^2 +\sum _{i=1}^n dx _i ^2,
\]
    where $H(u,x_1,\dots,x_n)= a(u) \sum _{i=1}^n x_i ^2 + \sum _{i=1}^n b_i (u) x_i +c(u)$,
    for some functions $a(u)$, $b_i(u)$, $c(u)$, and the function $f$ is a function of $u$ satisfying $f''(u)-\mu (f'(u))^2-n a(u)=0$.
    \end{enumerate}

\end{enumerate}
\end{theorem}

Depending on the character of $\nabla f$ we say that a quasi-Einstein manifold is {\it non isotropic} if $\|\nabla f\|\neq 0$ or {\it isotropic} if $\|\nabla f\|=0$. In the following sections we will study both cases separately.

The paper is organized as follows. In Section~\ref{se:2} we study some formulae and properties involving geometric objects of quasi-Einstein manifolds, specially the Ricci tensor, and introduce $pp$-waves. The analysis of case $(ii)(a)$ in Theorem~\ref{mainth} is carried out in Section~\ref{se:3}. In Section~\ref{se:4} we obtain some results on the isotropic case and, finally, we study locally conformally flat quasi-Einstein $pp$-waves in Section~\ref{se:5} to complete the proof of Theorem~\ref{mainth}.

\section{Preliminaries}\label{se:2}

Let $(M,g,f,\mu)$ be a Lorentzian quasi-Einstein manifold. As observed in the Introduction, if $\mu=-\frac 1n$, the manifold is globally conformally equivalent to a space form. In what follows we are going to introduce some results and definitions that we will use in the subsequent sections to prove Theorem \ref{mainth}-$(ii)$.

\subsection{General formulae.}
Although the next result was given in \cite{CMMR}, we include a sketch of a proof following a different strategy in order to make the paper as self-contained as possible.

\begin{lemma}
A Lorentzian quasi-Einstein manifold $(M,g,f,\mu)$ satisfies
\begin{align}
&\tau+\Delta f-\mu \|\nabla f\|^2 =(n+2)\lambda,\label{condquaeins1}\\
&\nabla \tau=2(\lambda-(n+2)\lambda \mu+\mu (1-\mu) \|\nabla f\|^2-\mu \tau)\nabla f+(\mu-1) \nabla \|\nabla f\|^2.\label{condquaeins2}
\end{align}
\end{lemma}
\proof
Equation~\eqref{condquaeins1} is obtained by simply contracting equation (\ref{q-einstein}).
Taking into account the contracted Second Bianchi identity $\nabla \tau=2 \text{div}(\rho)$ and the Bochner formula $\text{div}(\nabla \nabla f)=\rho(\nabla f)+\nabla \Delta f$ we compute the divergence of Equation~(\ref{q-einstein}):
\begin{align*}
0=&\text{div}(\lambda g)(X)\\
 =&\text{div}(\rho+\text{Hes}_f-\mu df\otimes df)(X)\\
 =&\frac 12 g(\nabla \tau,X)+\rho(\nabla f,X)+g(\nabla \Delta f,X)\\
 &-\mu g((\text{div}\nabla f)\nabla f,X)-\mu g(\nabla_{\nabla f}\nabla f,X).
\end{align*}

Covariantly differentiating Equation~\eqref{condquaeins1} we get $\nabla \tau=-\nabla \Delta f+\mu \nabla \|\nabla f\|^2$. Using this equation and that $\nabla_{\nabla f}\nabla f=\frac 12 \nabla \|\nabla f\|^2$ we see that
\begin{align*}
0=&\frac 12 g(\nabla \tau,X)+\lambda g(\nabla f, X)+\mu \|\nabla f\|^2g(\nabla f,X)-g(\nabla_{\nabla f} \nabla f,X)\\
 &+g(\mu \nabla \|\nabla f\|^2,X)-g(\nabla \tau,X)-\mu g(\Delta f \nabla f,X)-\mu g(\nabla_{\nabla f}\nabla f,X)\\
= & g\left(-\frac 12 \nabla \tau+(\lambda+\mu \|\nabla f\|^2-\mu \Delta f)\nabla f+\left(\frac{\mu-1}{2}\right)\nabla \|\nabla f\|^2,X\right).
\end{align*}

Now we replace $\Delta f$ by $(n+2)\lambda+\mu \|\nabla f\|^2-\tau$
to obtain Equation~\eqref{condquaeins2}.
\qed

\begin{remark}\label{nablatauisotropic}\rm
If the quasi-Einstein manifold $(M,g,f,\mu)$ is isotropic, i.e., $\|\nabla f\|=0$, then Equation~\eqref{condquaeins2} reduces to
\[
\nabla \tau=2(\lambda-\mu ((n+2)\lambda-\tau))\nabla f.
\]
Also note that, from Equation~\eqref{q-einstein}, one can write the Ricci operator ($g(\text{Ric} (X),Y)=\rho(X,Y)$) in the direction of $\nabla f$ as
\[
2Ric(\nabla f)=2\lambda \nabla f+2\mu \|\nabla f\|^2 \nabla f-\nabla \|\nabla f\|^2,
\]
so, if $\|\nabla f\|=0$, then $\text{Ric}(\nabla f)=\lambda \nabla f$ and $\nabla f$ is an eigenvector of the Ricci operator associated to the eigenvalue $\lambda$.
\end{remark}


\subsection{Some curvature properties of locally conformally flat quasi-Einstein manifolds.}
Let $(M,g)$ be an arbitrary $(n+2)$-dimensional Lorentzian manifold. If $n\geq 2$, local conformal flatness is characterized by the fact that the curvature tensor is given by
\begin{equation}\label{curv-general}
\begin{array}{rcl}
R(X,Y,Z,T)& =& \displaystyle\frac{\tau}{n(n+1)}\left\{g(X,T)g(Y,Z)-g(X,Z)g(Y,T)\right\}\\
\noalign{\medskip}
            &&\phantom{\frac{\tau}{n(n+1)}\left\{\right.}+\frac{1}{n}\left\{\rho(X,Z)g(Y,T)+\rho(Y,T)g(X,Z)\right.\\
\noalign{\medskip}
            &&\phantom{\frac{\tau}{n(n+1)}+\frac{1}{n}\left\{\right.}\left.-\rho(X,T)g(Y,Z)-\rho(Y,Z)g(X,T)\right\}.
\end{array}
\end{equation}
since the Weyl tensor is zero. If $n=1$, \eqref{curv-general} is always satisfied and does not characterize local conformal flatness. Let $C=\frac 1n\left(\rho-\frac{\tau}{2(n+1)g} \right)$ denote the Schouten tensor. For any $n\geq 1$, if $(M,g)$ is locally conformally flat then $C$ is a Codazzi tensor, i.e., its covariant derivative is totally symmetric and, moreover, this property characterizes local conformal flatness if $n=1$.

We proceed as in \cite{FG} to begin the study of the spectrum of the Ricci tensor.
\begin{lemma}\label{le:3}
Let $(M,g,f,\mu)$ be a locally conformally flat quasi-Einstein manifold of dimension $n+2$. Then, if $\mu\neq -\frac{1}{n}$, $\nabla f$ is an eigenvector of the Ricci operator.
\end{lemma}
\proof
Since $(M,g)$ is locally conformally flat the Schouten tensor is Codazzi, so \begin{equation}\label{codazzi1}
(\nabla_X \rho)(Y,Z)-\frac{X(\tau)}{2(n+1)}g(Y,Z)=(\nabla_Y \rho)(X,Z)-\frac{Y(\tau)}{2(n+1)}g(X,Z),
\end{equation}
for all vector fields $X,Y,Z$.


Using Equation~\eqref{q-einstein} we can write
\begin{align*}
(\nabla_X \rho)(Y,Z)=&-(\nabla_X Hes_f)(Y,Z)+\mu(\nabla_X df\otimes df)(Y,Z)\\
  =&-X(Hes_f)(Y,Z)+g(\nabla_{\nabla_XY}\nabla f,Z)+g(\nabla_Y\nabla f,\nabla_XZ)\\
   &+\mu\left(X(df\otimes df (Y,Z))-g(\nabla f,\nabla_XY)Z(f)-g(\nabla f,\nabla_XZ)Y(f)\right)\\
  =&-g(\nabla_X \nabla_Y \nabla f,Z)+g(\nabla_{\nabla_X Y}\nabla f,Z)\\
   &+\mu (\text{Hes}_f(X,Y)g(\nabla f,Z)+\text{Hes}_f(X,Z)g(\nabla f,Y)).
\end{align*}

We substitute this expression in \eqref{codazzi1} to get
\begin{align*}
&-g(\nabla_X \nabla_Y \nabla f,Z)+g(\nabla_{\nabla_X Y}\nabla f,Z)-\frac{X(\tau)}{2(n+1)}g(Y,Z)+\mu \text{Hes}_f(X,Z)g(\nabla f,Y)\\
=&-g(\nabla_Y \nabla_X \nabla f,Z)+g(\nabla_{\nabla_Y X}\nabla f,Z)-\frac{Y(\tau)}{2(n+1)}g(X,Z)+\mu \text{Hes}_f(Y,Z)g(\nabla f,X).
\end{align*}
Reorganizing the terms of this expression and using that the curvature tensor is given by $R(X,Y)Z=\nabla_X\nabla_YZ-\nabla_Y\nabla_XZ-\nabla_{[X,Y]}Z$, we obtain
\begin{align}\label{curv-simplify}
R(X,Y,Z,\nabla f)=&-\frac{X(\tau)}{2(n+1)}g(Y,Z)+\frac{Y(\tau)}{2(n+1)}g(X,Z)\\
     &+\mu (\text{Hes}_f(X,Z)g(\nabla f,Y)-\text{Hes}_f(Y,Z)g(\nabla f,X)).\nonumber
\end{align}

%

We choose vector fields $Z=\nabla f$ and $X$ such that $g(X,\nabla f)=1$ to get that
\[
0=\frac{Y(\tau)}{2(n+1)}-\mu \text{Hes}_f(Y,\nabla f)
\]
for all $Y\perp \nabla f$. Now, from Equation~\eqref{condquaeins2} we have that $Y(\tau)=2(\mu-1) \text{Hes}_f(Y,\nabla f)$ so
\[
0=-\frac{n\mu+1}{n+1}\text{Hes}_f(Y,\nabla f).
\]
Hence, either $\mu=-\frac{1}{n}$ or $\nabla f$  is an eigenvector of the Hessian operator $\text{hes}_f(X)=\nabla_X\nabla f$. Assume the latter, from \eqref{q-einstein} we have that $\rho(Y,\nabla f)=-\text{Hes}_f(Y,\nabla f)$ and therefore
\[
0=\frac{n\mu+1}{n+1}\rho(Y,\nabla f),
\]
showing that $\nabla f$ is also an eigenvector of the Ricci operator, unless $\mu=-\frac{1}{n}$.
\qed

\subsection{$pp$-waves.}
A Lorentzian manifold is said to be indecomposable if the metric is degenerate on any invariant proper subspace of the holonomy group. Indecomposable but not irreducible Lorentzian manifolds admit a parallel degenerate line field $\mathcal{D}$ and the curvature tensor satisfies
\[
R(\mathcal{D},\mathcal{D}^\perp,\cdot,\cdot)=0, \qquad R(\mathcal{D},\mathcal{D},\cdot,\cdot)=0, \qquad \text{and } \qquad R(\mathcal{D}^\perp,\mathcal{D}^\perp,\mathcal{D},\cdot)=0.
\]
Among manifolds with these properties, those with $\mathcal{D}$ spanned by a parallel null vector field and satisfying $R(\mathcal{D}^\perp,\mathcal{D}^\perp,\cdot,\cdot)=0$ constitute the family described in the following definition and play an important role both in Physics and Mathematics.
\begin{definition}\label{def:pp-wave}
A manifold $(M,g)$ locally isometric to $(\mathbb{R}^2\times \mathbb{R}^n,g_{ppw})$ with metric given in local coordinates $(u,v,x_1,\dots,x_{n})$ by
\begin{equation}\label{pp}
g_{ppw}=2dudv +H(u,x_1,\dots,x_{n}) du^2 +\sum _{i=1}^n dx _i ^2,
\end{equation}
where $H(u,x_1,\dots,x_n)$ is an arbitrary smooth function, is called a $pp$-wave. Moreover, if $H(u,\cdot)$ is a quadratic form in $\mathbb{R}^n$ then $(M,g)$ is said to be a {\it plane wave}.
\end{definition}
In Section~\ref{se:5} we will use the following result given in \cite{leistner} that characterizes $pp$-waves in terms of the Ricci tensor.
\begin{lemma}\label{lemma:leistner-ppwave}
Let $(M,g)$ be a Lorentzian manifold admitting a parallel null line field $\mathcal{D}$ so that $R(\mathcal{D}^\perp,\mathcal{D}^\perp,\cdot,\cdot)=0$. If the image of the Ricci operator $\text{Ric}$ is totally isotropic then $(M,g)$ is a $pp$-wave.
\end{lemma}

\section{Non isotropic locally conformally flat quasi-Einstein manifolds}\label{se:3}

In this section we show that in a neighborhood of any point where $\|\nabla f\|\neq 0$ the underlying manifold has the local structure of a warped product, thus proving Theorem \ref{mainth}-$(ii)-(a)$.

\begin{lemma}\label{lemma:non-isotropic}
Let $(M,g,f,\mu)$ be a locally conformally flat Lorentzian quasi-Einstein manifold with $\|\nabla f\|_P\neq 0$ in some point $P\in M$. Then, if $\mu\neq -\frac{1}{n}$,  $(M,g)$ is a warped product of a real interval and a space of constant sectional curvature on a neighborhood of $P$.
\end{lemma}
\proof
Since $\|\nabla f\|_P\neq 0$, $\|\nabla f\|\neq 0$ on a neighborhood $\mathcal{U}$ of $P$. Consider the unit vector $V=\frac{\nabla f}{\|\nabla f\|}$ on $\mathcal{U}$, which can be timelike or spacelike (we set $\varepsilon=g(V,V)=\pm 1$). Consider a local orthonormal frame $\{V=E_0,E_1,\dots,E_{n+1}\}$ and set $\varepsilon_i=g(E_i,E_i)$.

From Lemma~\ref{le:3} we have that $\rho(E_i,V)=\text{Hes}_f(E_i,V)=0$ for all $i=1,\dots,n+1$, hence from Equation~\eqref{condquaeins2} we obtain
\[
E_i(\tau)=2\|\nabla f\|(1-\mu)\rho(V,E_i)=0.
\]
We compute $R(V,E_i,E_i,V)$ in Equations~\eqref{curv-simplify} and \eqref{curv-general} to see that
\begin{equation*}
-\frac{V(\tau)}{2(n+1) \|\nabla f\|}\varepsilon_{i}-\mu \text{Hes}_f(E_i,E_i)\varepsilon=\frac{\tau}{n(n+1)}\varepsilon_{i}\varepsilon-\frac 1n \rho(V,V)\varepsilon_{i}-\frac 1n \rho(E_i,E_i)\varepsilon,
\end{equation*}
from where
\[
\left(\mu+\frac 1n\right)\text{Hes}_f(E_i,E_i)\varepsilon=\left(-\frac{\tau}{n(n+1)}\varepsilon+\frac 1n \rho(V,V)+\frac 1n\lambda \varepsilon-\frac{V(\tau)}{2(n+1) \|\nabla f\|}\right)\varepsilon_{i}.
\]
This shows that, when $\mu\neq -\frac{1}{n}$ the level sets of $f$ are totally umbilical hypersurfaces. Hence, as the normal foliations are totally geodesic ($g(\nabla_{\nabla f} \nabla f,E_i)=\text{Hes}_f(\nabla f,E_i)=0, \, \forall i=1,\dots,n+1$), $(M,g)$ decomposes locally as a twisted product (see  \cite{PR}). Now, since  $\rho(V,E_i)=0$ for all $i=1,\dots,n+1$,
the twisted product reduces to a warped product \cite{FGKU}. In conclusion $(M,g)$ is locally a warped product
$(I\times N,\varepsilon dt^2+\psi(t)^2 g_N)$, which is locally conformally flat by hypothesis and hence $(N,g_N)$ is a Riemannian or a Lorentzian manifold of constant sectional curvature~\cite{BGV}.
\qed

\begin{remark}\rm
Note that in Riemannian signature a quasi-Einstein manifold satisfies the hypothesis of Lemma~\ref{lemma:non-isotropic} and analogous arguments to those given here apply. Thus one has that locally conformally flat quasi-Einstein manifolds are either conformally equivalent to a space form or locally isometric to a warped product of a real interval and a space of constant sectional curvature. We refer to \cite{CMMR} for a different proof of this result.
\end{remark}

\section{Isotropic locally conformally flat quasi-Einstein manifolds}\label{se:4}

In order to prove Theorem~\ref{mainth}-$(ii)-(b)$, in this section we are going to study isotropic locally conformally flat Lorentzian quasi-Einstein manifolds.

\begin{lemma}\label{lemma:main-isotropic}
Let $(M,g,f,\mu)$ be an isotropic locally conformally flat Lorentzian quasi-Einstein manifold. Then, if $\mu\neq -\frac 1n$, around any regular point of $f$, the manifold $(M,g)$ is locally a $pp$-wave.
\end{lemma}

\proof
Let $(M,g,f,\mu)$ be quasi-Einstein with $\|\nabla f\|=0$. We will see that $\nabla f$ spans a parallel null line field, but we start by analyzing the Ricci tensor. First we choose an appropriate basis to work with. Set $V=\nabla f$. Since $V$ is a null vector, there exist orthogonal vectors $S$, $T$ satisfying $g(S,S)=-g(T,T)=\frac 12$ such that $V=S+T$. Define $U=S-T$, which is a null vector such that $g(U,V)=g(S,S)-g(T,T)=1$, and consider a basis $\{U,V,E_1,\dots,E_n\}$, with $g(E_i,E_i)=1$ for all $i=1,\dots,n$. We begin the study of the Ricci tensor by noting that, since $\operatorname{Ric}(V)=\lambda V$ by Remark~\ref{nablatauisotropic}, we have
\[
\rho(V,V)=0,\,\rho(U,V)=\lambda \text{ and } \rho(V,E_i)=0\quad  \forall i=1,\dots,n.
\]
We compute $R(U,E_i,E_j,V)$ both in Equation~\eqref{curv-simplify} and in Equation~\eqref{curv-general} to see that
\begin{equation}\label{eq:isotropic-relation}
\frac{-U(\tau)}{2(n+1)}\delta_{ij}-\mu \text{Hes}_f(E_i,E_j)
     = \frac{\tau}{n(n+1)}\delta_{ij}-\frac{1}{n}\lambda\delta_{ij}-\frac{1}{n}\rho(E_i,E_j).
\end{equation}
From Equation~\eqref{q-einstein} we get that $\text{Hes}_f(E_i,E_j)=\lambda\delta_{ij}-\rho(E_i,E_j)$ and from Remark~\ref{nablatauisotropic} that $U(\tau)=2(\lambda-\mu((n+2)\lambda-\tau))$. Hence, since $\mu\neq -\frac1n$, from Equation~\eqref{eq:isotropic-relation} we conclude that, if $i\neq j$, then $\rho(E_i,E_j)=0$, whereas $\rho(E_i,E_i)=\frac{\tau-\lambda}{n+1}$ for all $i=1,\dots, n$. Now, compute the scalar curvature
\[
\tau=2\rho(U,V)+n\rho(E_i,E_i)=2\lambda+n \frac{\tau-\lambda}{n+1},
\]
to see that
\[
\tau=(n+2)\lambda,
\]
so $\tau$ is constant. Moreover, since $\nabla \tau=0$, by Remark~\ref{nablatauisotropic} we have that $\lambda-\mu((n+2)\lambda-\tau)=0$ and then $\lambda=0=\tau$.
From \eqref{curv-general} we compute
\[
\rho(U,E_i)=R(U,V,E_i,V)+\sum_{j\neq i} R(U,E_j,E_i,E_j)=\frac{n-1}n \rho(U,E_i).
\]
Thus $\rho(U,E_i)=0$ for $i=1,\dots,n$ and the only possibly nonzero term of the Ricci tensor is $\rho(U,U)$. Therefore the Ricci operator is two-step nilpotent and its image is totally isotropic.

Now we are going to show that the line field $\mathcal{D}=\operatorname{span}\{V\}$ is parallel by checking that $V$ is a recurrent vector field. Thus, from Equation~\eqref{q-einstein} we compute
\begin{align*}
&\text{hes}_h(U)=-Ric(U)+\mu df\otimes df(U)=-\rho(U,U) U+\mu U,\\
&\text{hes}_h(V)=-Ric(V)+\mu df\otimes df(V)=0,\\
&\text{hes}_h(E_i)=-Ric(E_i)+\mu df\otimes df(E_i)=0,
\end{align*}
to see that $\nabla_X \nabla f=\omega(X)\nabla f$ for the $1$-form $\omega$ given by $\omega(U)=\mu-\rho(U,U)$, $\omega(V)=0$, $\omega(E_i)=0$ for all $i=1,\dots,n$.

It is now easy to see from \eqref{curv-general} and the information on the Ricci tensor that \[
R(\mathcal{D}^\perp,\mathcal{D}^\perp,\cdot,\cdot)=0.
\]
Hence, by Lemma~\ref{lemma:leistner-ppwave}, it follows that $(M,g)$ is a $pp$-wave.
\qed

\begin{remark}\rm
Note that although $(M,g)$ is a $pp$-wave, and hence it admits a null parallel vector field, $\nabla f$ is not in general parallel.
\end{remark}

\section{Locally conformally flat quasi-Einstein $pp$-waves}\label{se:5}

Let $(M,g_{ppw})$ be a $pp$-wave given in local coordinates as in \eqref{pp}. Denote $\partial_u=\frac{\partial}{\partial u}$, $\partial_v=\frac{\partial}{\partial v}$ and $\partial_i=\frac{\partial}{\partial x_i}$ for all $i=1,\dots,n$.
The possibly non-zero components of the curvature tensor are given, up to the usual symmetries, by
\begin{equation}\label{curvcomp}
R(\partial_u,\partial_i,\partial_u,\partial_j)=-\frac 12 \partial_{ij}^2 H,  \quad i,j=1,\dots,n.
\end{equation}
Hence the only possibly non-zero component of the Ricci tensor is
 \begin{equation}\label{rhoei}
\rho (\partial_u,\partial_u) = -\frac{1}{2} \sum_{i=1}^n \partial_{ii}^2 H.
\end{equation}
Thus, a $pp$-wave is Einstein, and moreover Ricci flat, if and only if the space-Laplacian of the defining function $H$ vanishes identically.
As a consequence of expressions \eqref{curvcomp} and \eqref{rhoei}, a $pp$-wave as in \eqref{pp} is locally conformally flat if
\begin{equation}\label{cflat}
H(u,x_1,\dots,x_n)= a(u) \sum _{i=1}^n x_i ^2 + \sum _{i=1}^n b_i (u) x_i +c(u),
\end{equation}
for arbitrary smooth functions $a,b_1,\dots,b_n,c$ on the variable $u$. Furthermore, note that
in this case the Ricci tensor is only a function of $u$, $\rho(\partial_u,\partial_u)=-n a(u)$. Recall that a $pp$-wave whose defining function $H(u,\cdot)$ is a quadratic form on the variables $x_1,\dots,x_n$ is called a plane wave. Thus, one can change variables in Equation~\eqref{cflat} to see that locally conformally flat $pp$-waves are a particular family of plane waves.

The following result characterizes locally conformally flat quasi-Einstein $pp$-waves.
\begin{lemma}\label{lemma:lcf-quasi-Einstein-pp-wave}
Let $(M,g_{ppw},f,\mu)$ be a non-trivial locally conformally flat quasi-Einstein $pp$-wave with $g_{ppw}$ given as in \eqref{pp}. Then one of the following holds:
\begin{enumerate}
\item[(i)] $\mu=0$ and $(M,g_{ppw},f)$ is a steady gradient Ricci soliton with $f$ a function of $u$ satisfying $f''(u)=-\rho(\partial_u,\partial_u)=n a(u)$.
\item[(ii)] $\mu=-\frac{1}{n}$ and $(M,g_{ppw})$ is conformally equivalent to a manifold of constant sectional curvature $c\leq 0$. Moreover, if $M$ is $3$-dimensional or $\|\nabla f\|=0$, it is conformally equivalent to a flat manifold.
\item[(iii)] $\mu\neq 0$ and $\mu\neq -\frac{1}{n}$, then the function $f$ is given by
\[
f(u,v,x_1,\dots,x_n)=-\frac{1}{\mu}\text{log}\{f_0(u)\}
\]
with
\[
f''_0(u)=\rho(\partial_u,\partial_u)\mu f_0(u)=-n a(u)\mu f_0(u).
\]
\end{enumerate}
\end{lemma}
\proof
If $\mu=0$ then $(M,g_{ppw},f)$ is a gradient Ricci soliton and the result follows from \cite{BGG}. If $\mu=-\frac{1}{n}$ then it follows
that $(M,g_{ppw})$ is conformally equivalent to an Einstein manifold and, since $(M,g_{ppw})$ is locally conformally flat, it is conformally equivalent to a manifold of constant sectional curvature. Furthermore we also have that $\tau=\lambda=0$ (see \cite{BGG1}) and from Equation~\eqref{condquaeins1} we obtain that $\Delta f=\frac{-1}{n}\|\nabla f\|^2$. Hence the conformal factor
reduces to $\frac1n\left(\frac1n-1\right)\|\nabla f\|^2$, which is zero whenever $n=1$ or $\|\nabla f\|=0$, or otherwise positive.

If $\mu\neq 0$ and $\mu\neq -\frac{1}{n}$, let $f(u,v,x_1,\dots,x_n)$ be an arbitrary function. The gradient of $f$ with respect to \eqref{pp} is given by
$\nabla f=(\partial_vf,\partial_u f-H \partial_vf,\partial_1f,\dots,\partial_nf)$
and thus \eqref{q-einstein} reduces to the following equations:
\begin{equation}\label{eq1}
\frac 12\underset{i=1}{\overset{n}{\sum }} \partial_iH\, \partial_if+\partial^2_{uu}{f}-\frac 12\partial_{u}H\, \partial_v f+\rho(\partial_u,\partial_u)-\mu(\partial_uf)^2=\lambda H,
\end{equation}
\begin{align}
&\partial^2_{ui}f-\frac 12\partial_iH\, \partial_v f-\mu\partial_if\partial_u f=0,& \qquad 1\leq i\leq n,\label{eq2}
\\
\noalign{\smallskip}
&\partial_{ii}^2f-\mu(\partial_if)^2=\lambda,& \qquad 1\leq i\leq n,\label{eq3}
\\
\noalign{\smallskip}
&\partial_{uv}^2f-\mu\partial_uf \, \partial_vf=\lambda,&\label{eq4}
\\
\noalign{\smallskip}
&\partial^2_{vv}f-\mu(\partial_vf)^2=0,&\label{eq5}
\\
\noalign{\smallskip}
&\partial^2_{vi}f-\mu\partial_if \, \partial_vf=0,& \qquad 1\leq i\leq n,\label{eq6}
\\
\noalign{\smallskip}
&\partial^2_{ij}f-\mu\partial_if \, \partial_jf=0,& \qquad 1\leq i\neq j\leq n.\label{eq7}
\end{align}


First we are going to show that $\partial_v f=0$ and that $\lambda=0$. We argue as follows. Differenciate \eqref{eq4} with respect to $v$ to see that $\partial_{uvv}^3\,f-\mu\,\partial^2_{uv}f \, \partial_{v}f-\mu\,\partial_{u}f  \partial^2_{vv}f=0$. Differentiate \eqref{eq5} with respect to $u$ to see that $\partial_{uvv}^3f=2\mu\, \partial_vf\partial^2_{uv}f$ and substitute the third order term to get that $\mu\, \partial_vf\partial^2_{uv}f-\mu\,\partial_{u}f \partial^2_{vv}f=0$. Now use \eqref{eq4} and \eqref{eq5} to reduce higher order terms to first order terms and get $\lambda\,\mu\,\partial_v f=0$. Since $\mu\neq 0$, either $\lambda=0$ or $\partial_v f=0$. If $\partial_v f=0$ then, as a consequence of \eqref{eq4}, we also have that $\lambda=0$. Hence assume $\lambda=0$. Differentiate \eqref{eq2} with respect to $x_i$ to get
\begin{equation}\label{eq:intermedia}
\partial^3_{uii} f-\frac12\left(\partial^2_{ii}H\,\partial_v f+\partial_iH\,\partial^2_{vi} f\right)-\mu\left(\partial^2_{ii}\,f\partial_u f+\partial_if\,\partial^2_{ui}f\right)=0.
\end{equation}
Differentiate \eqref{eq3} with respect to $u$ to see that $\partial^3_{uii} f=2\mu\,\partial_if\partial^2_{ui}f$. Simplify \eqref{eq:intermedia} using \eqref{eq2}, \eqref{eq3} and \eqref{eq6} to obtain
\[
\frac12 \partial^2_{ii}\,H\partial_v f=-\lambda\, \mu\, \partial_u f=0.
\]
Since $\partial_{ii}^2H\neq 0$ for some $i$ (otherwise the manifold is flat), we conclude that $\partial_v f=0$.

From equations $\partial_{ii}^2f-\mu(\partial_if)^2=0$ we get that
\[
f(u,v,x_1,\dots,x_n)=-\frac{1}{\mu}\text{log}\{f_0(u)+\sum_i \kappa_i x_i\},
\]
for some constants $\kappa_i$ and some function $f_0(u)$. Substitute $H$ by the expression in \eqref{cflat} and simplify to see that equations \eqref{eq1}-\eqref{eq7} reduce to the condition
\[
\mu\,\rho(\partial_u,\partial_u)\, \left(f_0(u)+\sum_i x_i \kappa_i\right)=f''_0(u)+a(u)\sum_i \kappa_i x_i+\frac12 \sum_i \kappa_i b_i(u).
\]
Since $\rho(\partial_u,\partial_u)=-n a(u)$, differentiate with respect to $x_i$ to see that
\[
-n\mu a(u) \kappa_i=a(u) \kappa_i.
\]
Since $a(u)\neq 0$ (otherwise the manifold is flat) and $\mu\neq -\frac{1}{n}$, we conclude that $\kappa_i=0$ for all $i=1,\dots,n$ and the result follows.
\qed

\medskip

{\it Proof of Theorem~\ref{mainth}.} The result follows from
Lemmas~\ref{lemma:non-isotropic}, \ref{lemma:main-isotropic} and \ref{lemma:lcf-quasi-Einstein-pp-wave}.

\begin{remark}\rm
Although Theorem~\ref{mainth} is local in nature, under certain conditions one can ensure that there exist global examples. Some of the manifolds described in Lemma~\ref{lemma:lcf-quasi-Einstein-pp-wave} can be defined globally: if $\mu=0$, $(M,g_{ppw},f)$ is a global steady gradient Ricci soliton and, if $\mu\neq 0,-\frac1n$, a sufficient condition for $f$ to be defined globally is that $a(u)\neq 0$ for all $u\in \mathbb{R}$. For example, Cahen-Wallach symmetric spaces are plane waves as in Definition~\ref{def:pp-wave} with $H(u,x_1,\dots,x_n)=\sum_{i=1}^n a_i x_i^2$ where $a_1,\dots,a_n$ are real numbers and hence they provide examples of isotropic geodesically complete quasi-Einstein manifolds (see \cite{BGG1}).

In contrast, plane waves with $H(u,x_1,\dots,x_n)=\sum_{i,j=1}^n (a_{ij}u+b_{ij})x_ix_j$
where $(a_{ij})$ is a diagonal matrix with the diagonal elements $a_{11}\leq\dots \leq a_{nn}$ non-zero real numbers and $(b_{ij})$ is an arbitrary symmetric matrix of real numbers, are examples of two-symmetric manifolds, this is, manifolds satisfying $\nabla^2 R=0$ but $\nabla R\neq 0$ \cite{AlGalaev}, \cite{BSS2}. Their Ricci tensor is given by $\rho(\partial_u,\partial_u)=-\sum_{i=1}^n(b_{ii}+u a_{ii})$, which shows that the Ricci vanishes if $u=-(\sum_{i=1}^n b_{ii})(\sum_{i=1}^n a_{ii})^{-1}$.
Hence, a Sturm-Liouville argument shows that any solution of the differential equation $f_0''(u)-\mu \rho(\partial_u,\partial_u) f_0(u)=0$ has infinitely many zeroes and thus they result in local but not global examples of quasi-Einstein manifolds \cite{BGG1}.
\end{remark}

\bibliography{mibiblio}
\bibliographystyle{amsplain}

\end{document}